\newcommand{\comment}[1]{}

\documentclass[11pt]{article}

\newenvironment{proof}{\noindent {\em Proof}.\ }{\hspace*{\fill}$\halmos$\medskip}
\newcommand{\notcontrol}[1]{}
\newcommand{\notacc}[1]{#1}
\newcommand{\foracc}[1]{}

\comment{uncomment FOR ACC:

\documentclass[letterpaper,10 pt,conference]{ieeeconf}
\IEEEoverridecommandlockouts
\overrideIEEEmargins
\newcommand{\notcontrol}[1]{}
\newcommand{\notacc}[1]{}
\newcommand{\foracc}[1]{#1}
}

\usepackage{amssymb,color}
\usepackage{graphicx}

\newcommand{\halmos}{\rule{1ex}{1.4ex}}

\newtheorem{theorem}{Theorem}
\newtheorem{itlemma}{Lemma}[section] 
\newtheorem{itproposition}[itlemma]{Proposition}
\newtheorem{itcorollary}[itlemma]{Corollary}
\newtheorem{itremark}[itlemma]{Remark}
\newtheorem{itdefinition}[itlemma]{Definition}
\newtheorem{itexercise}[itlemma]{Exercise}
\newtheorem{itexample}[itlemma]{Example}

\newenvironment{lemma}{\begin{itlemma}\rm}{\end{itlemma}} 
\newenvironment{remark}{\begin{itremark}\rm}{\end{itremark}} 
\newenvironment{corollary}{\begin{itcorollary}\rm}{\end{itcorollary}}
\newenvironment{proposition}{\begin{itproposition}\rm}{\end{itproposition}}
\newenvironment{definition}{\begin{itdefinition}\rm}{\end{itdefinition}}
\newenvironment{exercise}{\begin{itexercise}\rm}{\end{itexercise}}
\newenvironment{example}{\begin{itexample}\rm}{\end{itexample}}

\newcommand{\twoif}[4]{
\left\{ \begin{array}{ll}#1&#2\\#3&#4\end{array}\right.
}

\newcommand{\text}[1]{\hbox{\rm \ #1\ \/}}
\newcommand{\be}[1]{\begin{equation}\label{#1}}
\newcommand{\ee}{\end{equation}}
\newcommand{\bl}[1]{\begin{lemma}\label{#1}}
\newcommand{\ble}[1]{\begin{lemmaex}\label{#1}}
\newcommand{\br}[1]{\begin{remark}\label{#1}}
\newcommand{\bt}[1]{\begin{theorem}\label{#1}}
\newcommand{\bd}[1]{\begin{definition}\label{#1}}
\newcommand{\bp}[1]{\begin{proposition}\label{#1}}
\newcommand{\bc}[1]{\begin{corollary}\label{#1}}
\newcommand{\bfact}[1]{\begin{fact}\label{#1}}
\newcommand{\ber}[1]{\begin{exercise}\label{#1}}
\newcommand{\bex}[1]{\begin{example}\label{#1}}
\newcommand{\bem}[1]{\begin{example}\label{#1}}  
\newcommand{\ec}{\mybox\end{corollary}}
\newcommand{\efact}{\mybox\end{fact}}
\newcommand{\eer}{\mybox\end{exercise}}
\newcommand{\eex}{\mybox\end{example}}
\newcommand{\eem}{\mybox\end{example}}
\newcommand{\el}{\mybox\end{lemma}}
\newcommand{\ele}{\mybox\end{lemmaex}}
\newcommand{\er}{\mybox\end{remark}}
\newcommand{\et}{\qed\end{theorem}}
\newcommand{\ed}{\mybox\end{definition}}
\newcommand{\ep}{\mybox\end{proposition}}
\newcommand{\epr}{\end{proof}}
\newcommand{\bpr}{\begin{proof}}

\newcommand{\ecs}{\end{corollary}}
\newcommand{\eers}{\end{exercise}}
\newcommand{\eexs}{\end{example}}
\newcommand{\eems}{\end{example}}
\newcommand{\els}{\end{lemma}}
\newcommand{\eles}{\end{lemmaex}}
\newcommand{\ers}{\end{remark}}
\newcommand{\ets}{\end{theorem}}
\newcommand{\eds}{\end{definition}}
\newcommand{\eps}{\end{proposition}}

\newcommand{\qed}{\hfill \mbox{$\halmos$}} 
\newcommand{\mybox}{\hfill \mbox{$\Box$}} 

\newcommand{\rref}[1]{(\ref{#1})}

\newcommand{\R}{\mathbb R}

\newcommand{\Ii}{\cal I}

\newcommand{\pulseu}{{\bf u}_{\tau , \alpha }}
\newcommand{\pu}[1]{{\bf u}_{#1, \alpha }}
\newcommand{\sigmao}{\sigma ^o}
\newcommand{\Uu}{{\cal I}_{\tau , \alpha }}

\newcommand{\Vv}{{\cal V}}
\newcommand{\wh}[1]{\widehat{#1}}

\newcommand{\Aa}{{\cal A}}
\newcommand{\Ss}{{\cal S}}
\newcommand{\Gg}{{\cal G}}
\newcommand{\Bb}{{\cal B}}
\newcommand{\Cc}{{\cal C}}
\newcommand{\Dd}{{\cal D}}
\newcommand{\Ee}{{\cal E}}
\newcommand{\Ff}{{\cal F}}
\newcommand{\Hh}{{\cal H}}
\newcommand{\Mm}{{\cal M}}
\newcommand{\Xx}{{\cal X}}
\newcommand{\Yy}{{\cal Y}}

\newcommand{\pulseub}{{\bf u}_{\tau , \alpha , \beta }}

\newcommand{\sys}{\sigma }
\newcommand{\hsys}{\widehat{\sigma }}
\newcommand{\ben}{\begin{enumerate}}
\newcommand{\een}{\end{enumerate}}
\newcommand{\bi}{\begin{itemize}}
\newcommand{\ei}{\end{itemize}}
\newcommand{\beq}{\begin{eqnarray}}
\newcommand{\eeq}{\end{eqnarray}}
\newcommand{\beqn} {\begin{eqnarray*}}
\newcommand{\eeqn} {\end{eqnarray*}}

\newcommand{\sysI}{{\cal S}^I_n}
\newcommand{\sysII}{{\cal S}^{II}_n}

\newcommand{\ssim}{\sim} 
\newcommand{\ioeq}[1]{\,\raisebox{-1.5ex}{$\stackrel{\equiv }{{\scriptstyle{#1}}}$}\,}
\newcommand{\ioeQ}{\equiv }

\newcommand{\pic}[2]{\includegraphics[scale=#1]{#2}}

\begin{document}

\title{\LARGE\bf 
\foracc{Remarks on }
Input Classes for Identification of Bilinear Systems}

\notacc{\maketitle}

\begin{center}
Eduardo D. Sontag%
\footnote{Partially supported by NSF grants DMS-0504557 and DMS-0614371}\\
Dept. of Mathematics\\
Rutgers University\\
New Brunswick, NJ 08903, USA\\
{\tt sontag@math.rutgers.edu, 732.445.3072}
\end{center}\medskip

\begin{center}
Yuan Wang%
\footnote{Partially supported by NSF grant DMS-0072620 and Chinese
   National Natural Science Foundation grant 60228003}\\
Dept. of Mathematical Sciences\\
Florida Atlantic University\\
Boca Raton, FL 33431, USA\\
{\tt ywang@math.fau.edu, 561.297.3317}
\end{center}\medskip

\begin{center}
Alexandre Megretski\\
Dept. of Elec. Engr. \& Computer Science\\
Mass. Inst. of Technology\\
Cambridge, MA02139\\
{\tt ameg@mit.edu, 617.253.9828}
\end{center}\medskip

\foracc{\maketitle}

\notacc{\bigskip}

\begin{abstract}

  This paper asks what classes of input signals are sufficient in order to
  completely identify the input/output behavior of generic bilinear systems.
  The main results are that step inputs are not sufficient, nor are single
  pulses, but the family of all pulses (of a fixed amplitude but varying
  widths) do suffice for identification.

\end{abstract}

\notacc{\newpage}

\section{Introduction}

\notcontrol{Contents typeset like this are not intended for a pure control-theory
  audience.  They will be commented-out from a paper for ACC or for a control
  journal.  Please ignore --- they are here only so that I can remind myself
  later of points to be made when talking to biologists and others! }

\notcontrol{
Control systems theory concerns the study of open dynamical systems, which
process time-dependent input signals (stimuli, ligands, controls, forcing
functions, disturbances, test signals) into output signals (responses,
measurements, 
read-outs, reporters).  Such input/output (i/o) systems may be studied by
themselves, or as components (subsystems, modules) of larger systems.
}

\begin{figure}[ht]
\begin{center}
\notacc{\pic{0.5}{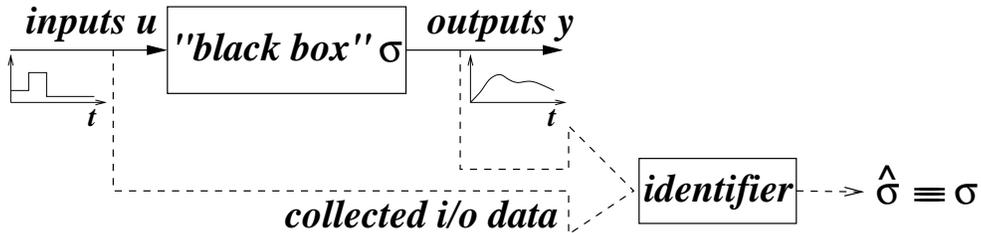}}
\foracc{\pic{0.3}{blackbox_id}}
\caption{identifying from i/o data}
\label{paradigm}
\end{center}
\end{figure}

In this paper, we address the following question: \emph{what types of input
signals are sufficient to completely identify the i/o behavior of a system?}
In other words, we look for classes ${\cal U}$ of inputs with the property that,
if a system $\sys$ is stimulated with the inputs from the set ${\cal U}$ and the
corresponding time record of outputs is recorded, then, on the basis
of the collected
information on inputs and outputs it is possible ---at least theoretically,
with no regard to computational effort, and in the absence of noise--- to
obtain
a system $\hsys$ which is equivalent to $\sys$ (Figure~\ref{paradigm}).
By ``equivalent'' we mean that the estimated system $\hsys$ will be
completely indistinguishable from the true system $\sys$ in its i/o behavior,
even when presented with inputs that do not belong to the restricted class
${\cal U}$ used for the identification
experiments.
Whether a certain class ${\cal U}$ of inputs is rich enough for identification is
heavily dependent, of course, on prior assumptions about the system $\sys$.

It is often very difficult to perform experiments in which arbitrary
input profiles are used.  Often, the only possible experiments are those in
which steps, i.e.\ constant inputs, are applied.  For example, in molecular
biology, a step input corresponds to subjecting a cell culture to a fixed
concentration of an extracelular ligand such as a drug or growth factor.
Sometimes somewhat more complex inputs, such as pulses (keep the input constant
at some level, then change it back to some default value) can be used, but
this is already not easy to implement, much less more complicated test signals.
This presents a theoretical challenge: how does one know if all possible
``identifiable'' information about the system can be obtained from such a
restricted class of experiments?  In the case of \emph{linear} models of
systems, this issue does not arise, because basically any single input, as
long as it is nonzero, for example a single step or a single pulse, suffices
for identification (or several inputs, one for each input channel, if the
system has multiple inputs; for simplicity we restrict here to single-input
systems). 
Note that we are thinking here of an experimental setup in which
observations are collected over time.  If, instead, only steady-state behavior
was observed, and not transients, then one input is not enough, even for
single-input linear systems.  In that case, one has to use multiple inputs,
such as steady-state measurements of responses to periodic signals at different
frequencies.

\notcontrol{
Linear theory is satisfactory for the analysis of ``weakly activated''
systems, that is to say the responses of systems to small signals and in
regimes of operation that are close to steady states.  It does not help
when dealing with systems in which nonlinear effects cannot be disregarded.}

For nonlinear systems, it is thus an interesting question whether constant
inputs or pulse inputs, or simple combinations of these input classes, suffice
for identification, as they do for linear systems.
In this paper, we show that for a large and interesting class of
nonlinear systems, that of \emph{bilinear systems}, constant inputs do not
suffice, but pulses do.

Bilinear systems constitute an appealing class of nonlinear
systems~\cite{bruni_et_al,mohler_book,mohler_paper,Elliott99}. 
While for linear systems the evolution of the states is only allowed to
depend on linear functions of the state variables and inputs, in bilinear
systems one also allows a linear dependence on products between input
and state variables. 
Bilinear systems can be easily described in linear-algebraic terms, and a
theory, in many respects analogous to the linear theory, can be developed for
the analysis of their input/output properties.
On the other hand, bilinear systems are theoretically
capable of approximating arbitrary
input/output behaviors on finite time-intervals~\cite{sussmann,fliess,rugh}.
They have been used to model chemical processes, electrical networks, power
plants, nuclear reactors, robotic manipulators, and many other systems in
engineering, chemistry, biology, economics, and other
fields~\cite{mohler_book}.
They can also be employed in order to model and analyze certain simple
enzymatic signaling cascades, when substrates are not too close to saturation
and thus Michaelis-Menten kinetics can be replaced by bilinear
expressions~\cite{heinrich02,chaves_sontag06}. 

Informally (see next section for precise definitions and statements), the main
results that we prove are as follows.  On the negative side, we show that step
inputs are not enough for identifying bilinear systems, nor do single pulses
suffice.  On the positive side, we show that the family of all pulses (of a
fixed amplitude but varying widths) do suffice for identification.  
To be precise, one must impose certain non-degeneracy conditions on the
classes of systems being considered, and, for the negative result, one wants to
avoid trivial counter-examples in order to say something interesting.  Thus,
all results are stated for generic classes of systems.

\notacc{%
Our techniques are based on realization theory.
\notcontrol{
Two extreme approaches to modeling systems (or particular subsystems) are as
follows~\cite{mct}:
\bi
\item
The \emph{black-box}, \emph{purely input/output}, or
\emph{phenomenological} description ignores all mechanistic detail and
characterizes behavior solely in terms of behavioral or stimulus-response
data.  This 
data may be summarized for example by transfer functions (for linear systems
and certain types of bilinear systems), Volterra, generating series, or
Wiener expansions (for general nonlinear systems), statistical models
that correlate inputs and outputs, neural network and other approximate models,
or machine learning prediction formalisms.
\item
The \emph{state-space}, \emph{internal}, or \emph{mechanistic} description,
in which all relevant state variables (such as, in a biological application,
proteins, mRNA, and metabolites) are specified, and evolution equations are
completely specified (in a biological application, this would mean the forms
for all reactions and reaction constants).
\ei
In between these two extremes, one finds combinations (``gray box'' view).
Available mechanistic information is incorporated into the model, and this is
supplemented by i/o data, which serves to impose constraints on (or, in a
Bayesian approach, to provide prior distributions of) unknown parameters to
be further identified.
In engineering applications, black-box modeling is common, for
example, in chemical process control, while state-space modeling is
often used in mechanical and aerospace engineering.
Many biological models involve combinations of both approaches.
Passing from one sort of representation to the other may be highly
nontrivial in practice, since data are incomplete and noisy.
{}From a theoretical point of view, however, the alternative methods of
description are closely related, in the sense that one may in principle go
from an internal to an external description, and, conversely, from the
knowledge of black-box, input/output i/o behavior, it is possible, at least in
principle, to characterize the set of all possible internal descriptions that
account for it.  This is the focus of realization theory (or ``reverse
engineering'').}
We make heavy use of ideas originally developed by Kalman for realizations by
linear systems, and refined by Isidori and later Fliess for bilinear systems.
}

The organization of this paper is as follows.  Section~\ref{sec-prelim}
gives the basic definitions, and reviews the background from bilinear
realization theory.  
We provide a self-contained discussion because, even though
the results proved in that section are not new, it is
hard to find references presented as needed here.
The main results are stated in Section~\ref{sec-statements}.
The proofs of the
negative results are given in Section~\ref{proof-neg-results} and the
proofs of the positive results are given in Section~\ref{proof-pos-results}.
The latter are inspired by Juang's very nice
paper~\cite{bilinear} (we thank David Angeli for bringing this paper to our
attention).
Although the bilinear identification algorithm given in~\cite{bilinear}
involves some ambiguous steps, such as taking (non-unique) logarithms of
matrices, we were able to adapt many of its basic ideas\notacc{;
we discuss in Section~\ref{sec-juang} connections to that work}.
Conclusions and some remarks are presented in Section~\ref{sec-final}.
\foracc{Many details have been left out due to space limitations.}

\section{Preliminaries}
\label{sec-prelim} 

\subsubsection*{Systems and I/O Equivalence}

\notacc{%
Our results will be for bilinear systems, which are defined by affine vector
fields on $\R^n$ and hence are described by matrices, but an abstract setup
will allow us to discuss some preliminary facts in more generality.
}
We consider single-input single-output initialized systems $\sys$,
in the usual sense of control theory~\cite{mct}:
\be{gensys}
\dot x(t) = f_0(x(t)) + u(t)f_1(x(t)),\;x(0)=x_0,\,\;y(t)=h(x(t))
\foracc{.}
\ee
\notacc{(we will drop the arguments ``$(t)$'' if clear from the context),
where $f_0$ and $f_1$ are smooth vector fields on a manifold $M$ and $h$ is an
output function $M\rightarrow \R$.
Inputs can be taken to be any Lebesgue-measurable, essentially
bounded, functions $u:[0,T_u]\rightarrow \R$, but there will be no difference in results
if one restricts to, for instance, piecewise continuous inputs.
We let $\Omega $ be the set of all inputs.
In principle, solutions are unique but only defined on some maximal interval
(which depends on the initial condition and the input), but for simplicity, and
since it is the case anyway for bilinear systems, we assume that solutions are
defined for all times $t\in [0,T_u]$ (forward completeness).
We let $\varphi(t,u)$, or just $x(t)$ if the input is clear from the context,
be the solution of~(\ref{gensys}) at time $t$, and $y(t)=h(\varphi(t,u))$ the
corresponding output.
When more than one system is being considered, we use appropriate notations;
for example, a system $\hsys$ will be described by equations
$\dot x = \wh{f}_0(x) + u\wh{f}_1(x)$,
$x(0)=\wh{x_0}$,
$y=\wh{h}(x)$.
}

Given two systems $\sys,\hsys$, and an input $u$, we say that 
$\sys,\hsys$ are \emph{indistinguishable} under $u$ if
$h(\varphi(t,u))=\widehat{h}(\widehat\varphi(t,u))$
for all $t\in [0,T_u]$.
\notacc{%
If $h(\varphi(t,u))\not= \widehat{h}(\widehat\varphi(t,u))$ for some $t$, we say that
$u$ distinguishes among the two systems $\sys,\hsys$;
in other words, the ``input/output experiment'' consisting of
perturbing the system with this input $u$ results in a different time-varying
observation for $\sys$ than for $\hsys$.

Given a subset ${\cal U}\subseteq \Omega $ of inputs, we say that the
two systems $\sys,\hsys$ are
}
\foracc{
Given a subset ${\cal U}\subseteq \Omega $ of inputs (where $\Omega $
denotes the set of all measurable and essentially bounded functions),
we say that the
two systems $\sys,\hsys$ are
}
\emph{input/output (i/o) equivalent with respect to all inputs in ${\cal U}$}
if no input in ${\cal U}$ is able to distinguish between the two systems
$\sys,\hsys$, i.e., provided that $h(\varphi(t,u))=\wh{h}(\wh{\varphi}(t,u))$
for each $u\in {\cal U}$ and for each $t\in [0,T_u]$.
We write $\sys\ioeq{{\cal U}}\hsys$ in that case.
\notacc{

}
In the special case in which ${\cal U}=\Omega $, we write simply $\sys\ioeQ\hsys$
and simply say that the two systems are i/o equivalent.
\notacc{%
That is to say, they cannot be
distinguished in any way whatsoever based on their ``back box'' input/output
behavior.
}

Let $\Sigma $ be a class of systems.
A subset ${\cal U}\subseteq \Omega $ of inputs is said to be 
\emph{sufficient for identifying systems in the class $\Sigma $}
if, for any pair of systems 
 $\sys,\hsys$ in $\Sigma $,
\[
\sys\ioeq{{\cal U}}\hsys \;\Rightarrow \; \sys\ioeQ\hsys \,.
\]
In other words: whenever $\sys$ is not i/o equivalent to $\hsys$, there must
exist some input in the set ${\cal U}$ which distinguishes among the two systems
$\sys,\hsys$.

\notacc{%
Linear systems (finite-dimensional, continuous-time) are those for which
$f_0$ is linear, $f_1$ is constant, and $h$ is linear, i.e.\ systems described
by equations
\be{e-linears}
\dot x = Ax+bu\,,\;\;x(0)=0\,,\;\;y=cx
\ee
with $A\in \R^{n\times n}$, $b\in \R^{n\times 1}$, and $c\in\R^{1\times n}$.
We often refer interchangeably to a linear system or its corresponding
triple of matrices $(A,b,c)$.
Linear systems can be identified by any single nonzero input on a nontrivial
interval, such as a constant function (a step), or a pulse.
}

\notcontrol{A sketch of proof is as follows.
We view the i/o behavior of a linear system in the frequency domain:
a multiplication $\hat y(s)=W(s)\hat u(s)$, where hats indicate Laplace
(shifted Fourier) transforms and $W(s)$ is the transfer function of the system
$\sys$.  Given \emph{any} nonzero input $u(t)$, its Laplace transform is again
nonzero, so one can simply estimate $W(s)$ by taking the quotient $\frac{\hat
y(s)}{\hat u(s)}$ of the Laplace transforms of the output and input.  Once
that $W(s)$ is so computed, the response to every other possible input $v(t)$
can be calculated as follows: (i) compute its Laplace transform, (ii) take the
product $\hat w(s) = W(s)\hat v(s)$, and finally (iii) use Laplace inversion
in order to obtain $w(t)$.
}

\subsubsection*{4-Tuples and Bilinear Systems}

We consider two slightly different classes of bilinear systems.
To define these classes, we first introduce \emph{4-tuples} as follows:
\be{4-tuple}
(A,N,b,c)\,\;\;\mbox{where}\;\;
A, N\in\R^{n\times n}, \; b\in\R^{n\times 1}, \;c\in\R^{1\times n}
\ee
(the integer $n$ is called the \emph{dimension} of the 4-tuple).

We say that a system~(\ref{gensys}) is a
\emph{bilinear system of type I} if $f_0$ is linear, $f_1$ is affine,
$x_0=0$, and $h$ is linear.  In other words, the system equations are:
\be{e-syso}
\begin{array}{clll}
\dot x &=& (A+uN)x + bu, & x(0) = 0,\\
y&=& cx,&{}
\end{array}
\ee
where $(A,N,b,c)$ is some 4-tuple as in~(\ref{4-tuple}).
We use a notation such as ``$\sigma ^o$'' to refer to a system of type I.
With some abuse of terminology, we also simply write
$\sigma ^o=(A, N, b, c)$.
Note that linear systems~(\ref{e-linears}) constitute the subclass of
bilinear systems of type I for which $N=0$.

We say that a system~(\ref{gensys}) is a
\emph{bilinear system of type II} if $f_0$ and $f_2$ are both
linear and $h$ is linear (but the initial state may be nonzero).
In other words, the system equations are:
\be{e-sys}
\begin{array}{clll}
\dot x &=& (A+uN)x, & x(0) = b,\\
y&=& cx,&{}
\end{array}
\ee
where $(A,N,b,c)$ is a 4-tuple as in~(\ref{4-tuple}).
Once again, we do not differentiate between a system of type II and its
associated 4-tuple when the meaning is clear from the context.

\subsubsection*{Similarity}
We say that two 4-tuples $(A, N, b, c)$ and $(\wh A, \wh N, \wh b, 
\wh c)$ are \emph{similar} if they have the same dimension $n$ 
and there exists an
invertible $n\times n$ matrix $T$ such that the following equalities hold:
\be{e-equiv}
A = T\wh AT^{-1},\;\;
N = T\wh NT^{-1},\;\;
b = T\wh b,\;\;
c = \wh cT^{-1}\,.
\ee
Note that, for linear systems
\foracc{
\be{e-linears}
\dot x = Ax+bu\,,\;\;x(0)=0\,,\;\;y=cx
\ee
}
($N=\wh N= 0_{n\times n}$), this reduces to the
familiar equivalence relation in linear systems theory.

We say that two bilinear systems $\sys, \wh{\sys}$ (both of type I or both of
type II) are \emph{similar} (or ``internally equivalent''), and we write
\[
\sys\ssim\wh{\sys}
\]
 if there is a change of variables $x=Tz$ such that
the equations of $\sys$ get transformed into those of $\wh\sys$.
For systems of type I, this means that
\[
(A+uN)x + bu \;=\;\dot x\;=\;T\dot z\;=\;
T[(\wh A+u\wh N)Tx + T\wh bu
\]
for all $x$ and $u$, and also $cx=cTz=\wh c z$; thus,
$\sys\ssim\wh \sys$ is the same as saying that
the 4-tuples 
$(A, N, b, c)$ and $(\wh A, \wh N, \wh b, \wh c)$ are similar.
An analogous statement holds for systems of type II.

An easy calculation shows that $\sys\ssim\wh \sys\Rightarrow 
\sys\ioeQ\wh \sys$, and 
a converse holds as well, under certain minimality assumptions, as discussed
below.

\notacc{%
\subsubsection*{Checking I/O Equivalence}

For analytic systems, input/output equivalence can be verified by checking
certain algebraic equalities, and there is no need to test all possible
inputs, as we discuss next.

For any smooth vector field $f$ on $\R^n$, and any smooth function
$h:\R^n\rightarrow \R$, the Lie derivative $L_fh:\R^n\rightarrow \R$ is defined as the function 
$(L_fh)(x) = \nabla h(x)\,f(x)$, where $\nabla h$ is the gradient of $h$.
(In differential-geometric terms, $L_fh$ is simply the value of the vector
field on $h$, when vector fields are viewed as derivations on spaces of
smooth functions.)
More generally, if $f_1,\ldots ,f_k$ are vector fields, the iterated derivative
$L_{f_k}\ldots L_{f_1}h$ is defined recursively by the formula
$L_{f_k}\left(L_{f_{k-1}}\ldots L_{f_1}h\right)$.

Suppose that $\sys$ and $\hsys$ are two systems~(\ref{gensys}) for which
the vector fields $f_0$ and $f_1$ are analytic and the function $h$ is
also analytic.  Then, $\sys\ioeQ\hsys$ if and only if
\be{equal-lie}
\left(L_{f_{i_k}}\ldots L_{f_{i_1}}h\right) (x_0)
\;=\;
\left(L_{\wh{f}_{i_k}}\ldots L_{\wh{f}_{i_1}}\wh{h}\right) (\wh{x_0})
\ee
for all sequences $(i_1,\ldots ,i_k)\in \{0,1\}^k$ and all $k\geq 0$.
(When $k=0$, (\ref{equal-lie}) says that $h(x_0)=\wh{h}(\wh{x_0})$.)
This is true because the expressions in~(\ref{equal-lie}) are the
coefficients of the Fliess generating series of the input/output
behavior associated to the respective systems, and the i/o behavior is in
one-to-one correspondence with the coefficients of the series,
see~\cite{wang_sontag_scl89}, Lemma 2.1.

For bilinear systems of types I or II, i/o equivalence amounts to
an equality of vectors.  
Indeed, take first systems of type I.
In this case, $f_0(x)=Ax$, $f_1(x)=Nx+b$, and $h(x)=cx$.
Therefore, one can see inductively that:
\[
\left(L_{f_{i_k}}\ldots L_{f_{i_1}}h\right)(x) \;=\;
\twoif
{\mbox{if}\; i_k=0:}{cA_{i_1}\ldots A_{i_k}x}
{\mbox{if}\; i_k=1:}{cA_{i_1}\ldots A_{i_k}x + cA_{i_1}\ldots A_{i_{k-1}}b}
\]
where $A_0=A$ and $A_1=N$.
In particular, for $x=x_0=0$, we have that 
$(L_{f_{i_k}}\ldots L_{f_{i_1}}h)(x)=0$
for all sequences with $i_k=0$, and
$(L_{f_{i_k}}\ldots L_{f_{i_1}}h)(x)=cA_{i_1}\ldots A_{i_{k-1}}b$
for all sequences with $i_k=1$.
}

Generally, given two 4-tuples $(A,N,b,c)$ and $(\wh{A},\wh{N},\wh{b},\wh{c})$,
let us say that they are i/o equivalent if
\be{io-bilin}
cA_{i_1}\ldots A_{i_{k}}b
\;=\;
\wh{c}\wh{A}_{i_1}\ldots \wh{A}_{i_{k}}\wh{b}
\ee
for all sequences of matrices $A_j$ picked out of $A$ and $N$,
including the ``empty'' sequence ($c b =\wh{c}\wh{b}$).
(It suffices to check sequences of length $n+\wh{n}$,
where $n,\wh{n}$ are the respective state-space dimensions;
cf.~\cite{brockett72,dalessandro_grasseli_isidori_siam73,isidori_TAC83,fliess73,springer89,cas89}.)

\notacc{Then the preceding discussion proves:}

\bl{ioeq-tuples-same}
Two systems $\sys$ and $\hsys$ of type I are i/o equivalent if
and only if the corresponding 4-tuples are i/o equivalent.
\el

For bilinear systems of type II, the same conclusion holds, in this
case because 
\[
\left(L_{f_{i_k}}\ldots L_{f_{i_1}}h\right)(b)
\;=\;
cA_{i_1}\ldots A_{i_k}b \;.
\]

\bl{ioeq-tuples-same-II}
Two systems $\sys$ and $\hsys$ of type II are i/o equivalent if
and only if the corresponding 4-tuples are i/o equivalent.
\el

\subsubsection*{Canonical Systems and Uniqueness}

A 4-tuple $(A,N,b,c)$ as in~(\ref{4-tuple})
will be said to be \emph{canonical}
provided that the following two properties hold:
\ben
\item
There is no proper subspace of $\R^n$ that contains $b$ and is invariant under
$x\mapsto Ax$ and $x\mapsto Nx$.
\item
There is no nonzero subspace of $\R^n$ that is contained in the
nullspace of $x\mapsto cx$ and is invariant under $x\mapsto Ax$ and $x\mapsto Nx$.
\een
\notacc{%
The first property can be equivalently expressed by saying that the set of
vectors of the form 
\be{span-reach}
A_{i_1}\ldots A_{i_k}b \,,
\ee
ranging over all matrix products with $A_j\in \{A,N\}$ (including $k=0$, i.e.,
$b$), or equivalently over all products of length $k$ at most $n-1$, must span
all of $\R^n$.
This property is often called ``span-reachability'' because, for bilinear systems,
it amounts to the requirement that the set of states reachable from the origin
span all of the state-space.  Similarly, the second property can be
equivalently expressed by the dual property that the span of the vectors
\[
A_{i_1}'\ldots A_{i_k}'c'
\]
(prime indicates transpose) be all of $\R^n$ (once again, length $\leq n-1$
suffices), and is an observability property for bilinear systems.
Canonical 4-tuples are also called ``minimal'' 
because~\cite{brockett72,dalessandro_grasseli_isidori_siam73,isidori_TAC83,fliess73,springer89,cas89}
they have minimal dimension among all other 4-tuples which are i/o equivalent
in the sense of~(\ref{io-bilin}); moreover, if a 4-tuple $(A,N,b,c)$
is not canonical, then~\cite{brockett72,dalessandro_grasseli_isidori_siam73,isidori_TAC83,fliess73,springer89,cas89} there is some 4-tuple
$(\wh{A},\wh{N},\wh{b},\wh{c})$ 
which is canonical and is so that~(\ref{io-bilin}) holds.
(We do not need in this paper the interpretations in terms of reachability
and observability, nor the minimality result.)

}
We will call a bilinear system $\sys$ (of type I or II) canonical if
the corresponding 4-tuple is canonical.

A very special case is that of linear systems~(\ref{e-linears}), i.e.\ systems
of type I with $N=0$).
Such a system $\sigma =(A,b,c)$ is canonical if and only if it is reachable and
observable in the usual sense of control theory~\cite{mct}.
The controllability matrix ${\cal R}(A, b)$ and the observability matrix ${\cal O}(A, c)$
are defined respectively by:
\beqn
{\cal R}(A, b) &=& \pmatrix{b& Ab& \ldots  & A^{n-1}b},\\
{\cal O}(A, c) &=& {\cal R}(A',c')' \;=\; \pmatrix{c'& A'c'& \ldots  & (A')^{n-1}c'}'
\eeqn
(prime indicates matrix transpose).
The system $\sigma $ is canonical iff both matrices have full rank $n$.

\subsubsection*{Similarity and I/O Equivalence}

We already remarked that 
 $\sys\ssim\wh\sys\Rightarrow \sys\ioeQ\wh\sys$
for any two bilinear systems (both of the same type).
Conversely, if both systems $\sys$ and $\wh\sys$ are canonical,
$\sys\ioeQ\wh\sys\Rightarrow \sys\ssim\wh\sys$.
Thus:
\be{canon-equiv}
\mbox{if $\sys$ and $\wh\sys$ are canonical,}\quad
\sys\ioeQ\wh\sys\;\Longleftrightarrow\;\sys\ssim\wh\sys \,.
\ee
This is a standard fact about bilinear
systems~\cite{brockett72,dalessandro_grasseli_isidori_siam73,isidori_TAC83,fliess73,springer89,cas89} (strictly speaking, these
references deal with discrete-time systems such as $x(t+1)=(A + u(t)N)x(t)$,
but the algebraic statement about 4-tuples is the same as in the
continuous-time case).  The proof is, in fact, completely analogous to the
proof for linear systems~\cite{mct}.
\notacc{For completeness, we provide a proof here:}

\bl{main-theo-canonical}
Suppose that the two 4-tuples $(A,N,b,c)$ and $(\wh{A},\wh{N},\wh{b},\wh{c})$
are canonical and i/o equivalent.  Then they are similar.
Moreover, the similarity transformation $T$ in~(\ref{e-equiv}) is unique.
\els

\notacc{%
\bpr
Pick any $\wh{x}\in \R^n$.
By the span-reachability property~(\ref{span-reach})
for $(\wh{A},\wh{N},\wh{b},\wh{c})$, there are real numbers
$\lambda _\alpha $, where $\alpha $ denotes sequences $(i_1,\ldots ,i_k)$ of length at most
$n-1$ (including the ``empty'' sequence) such that 
$\wh{x}=\sum_{\alpha }\lambda _\alpha \wh{A}_{\alpha }\wh{b}$,
where we denote $A_\alpha  = A_{i_1}\ldots A_{i_k}$ for $\alpha =(i_1,\ldots ,i_k)$.
Now define $T\wh{x}:=\sum_{\alpha }\lambda _\alpha A_{\alpha }b$.

There are many possible representations of a vector $\wh{x}$ as a linear
combination of the spanning set in~(\ref{span-reach})
for $(\wh{A},\wh{N},\wh{b},\wh{c})$, so to see that $T$ is
well-defined as a mapping we need to verify that 
if 
$\sum_{\alpha }\lambda _\alpha \wh{A}_{\alpha }\wh{b}=\sum_{\alpha }\lambda _\alpha '\wh{A}_{\alpha }\wh{b}$
then $\sum_{\alpha }\lambda _\alpha A_{\alpha }b=\sum_{\alpha }\lambda _\alpha 'A_{\alpha }b$.  By linearity,
it is enough to show that $\sum_{\alpha }\lambda _\alpha \wh{A}_{\alpha }\wh{b}=0$ $\Rightarrow $
$\sum_{\alpha }\lambda _\alpha A_{\alpha }b=0$.
Suppose that $\sum_{\alpha }\lambda _\alpha \wh{A}_{\alpha }\wh{b}=0$.
Then also $\wh{c}\wh{A}_\beta \sum_{\alpha }\lambda _\alpha \wh{A}_{\alpha }\wh{b}=0$ for any other
index 
$\beta $, or equivalently $\sum_{\alpha }\lambda _\alpha  \wh{c}\wh{A}_{\beta \alpha }b=0$,
where $\beta \alpha $ is the concatenation of the sequences $\beta $ and $\alpha $.
Now, i/o equivalence of the two given 4-tuples implies that
$\wh{c}\wh{A}_{\beta \alpha }\wh{b}=cA_{\beta \alpha }b$ for all indices,
and so also
$cA_{\beta }x=\sum_{\alpha }\lambda _\alpha  cA_{\beta \alpha }b=0$.
This holds for any index $\beta $, so, using the observability of
$(A,N,b,c)$, we conclude that $x=0$, as desired.
The mapping $T$ is obviously linear (by definition), and it is onto because of
the reachability of $(A,N,b,c)$, which means that every
$x\in \R^n$ can be written as $\sum_{\alpha }\lambda _\alpha A_{\alpha }b$ for some $\lambda _\alpha $'s.
To prove that $T$ is one-to-one, we simply reverse the argument used to prove
that $T$ was well-defined.
Uniqueness follows by the same argument.
\epr

By picking among all the possible linear combinations the one whose
coefficients have minimal Euclidean norm, one obtains
an explicit expression for $T$:
\[
T\;=\; {{\cal R}}\,\wh{{\cal R}}^{\#} \,,
\]
where $\#$ denotes matrix pseudoinversion,
${{\cal R}}$ is a matrix listing the products in~(\ref{span-reach})
of length $\leq n-1$, and $\wh{{\cal R}}$ lists the vectors in the same order
for the second 4-tuple.
}
For linear systems $\sigma =(A,b,c)$, the equivalence becomes
\be{lin-T}
T\;=\; {\cal R}(A,b)\, {\cal R}(\wh{A},\wh{b})^{-1}\,,
\ee
where $R$ is the usual reachability matrix (\cite{mct}, Theorem 27).

\subsubsection*{Generic Sets of Systems}

We will make statements about ``generic'' classes of systems, so we must
define this term carefully.  Genericity can be defined in many ways, for
example in probabilistic terms (a set is generic if it has ``probability one'')
or, as usual in mathematics, in terms of open dense sets.
In order to provide the strongest possible results, we combine both definitions
and say here that a subset $S$ of an Euclidean space $\R^\ell$ is
\emph{generic} provided that:
\bi
\item
the set $S$ has full measure, that is, the complement $S^c$ has Lebesgue
measure zero, and
\item
the set $S$ is open (and dense) in $\R^\ell$.
\ei
When dealing with sets of 4-tuples~(\ref{4-tuple}), we view such sets as
subsets of $\R^\ell$ with $\ell=2n^2+2n$.

When talking about genericity of classes of systems of type I or II,
we mean genericity of the sets of associated 4-tuples.
Specifically, if we let $\sysI$ be the class of $n$-dimensional bilinear
systems of type I, then we think of $\sysI$ as $\R^{2n^2+2n}$,
and similarly for the class $\sysII$ of 
$n$-dimensional bilinear systems of type II.

\section{Statements of Main Results}
\label{sec-statements}

For any $\alpha \in\R$ and any $\tau \ge 0$, let $\Uu$ denote
the class of all functions of the form
\[
\pulseub(t) 
= \left\{
\begin{array}{ll}
\alpha & \mbox{for $0\le t < \tau $,}\\
\beta  & \mbox{for $t\ge \tau $},
\end{array}\right.
\]
where $\beta $ is a constant. Let ${\bf u}_{\tau ,\alpha }$ 
denote the particular pulse function in $\Uu$ for which $\beta =0$, that
is, 
\be{e-pulseu}
\pulseu(t) = \left\{
\begin{array}{ll}
\alpha & \mbox{for $0\le t < \tau $,}\\
0 & \mbox{for $t\ge \tau $}.
\end{array}\right.
\ee
Note that in the special case when $\tau =0$, $\Ii_{\tau , \alpha }$
becomes the class of constant functions.
(There is a small ambiguity in that we have not specified the domain of
the inputs.  We can view these inputs as defined on some interval
$[0,T]$ with $T>\tau $; any such $T$ will give the same results.)

Let us now state the negative main results of this paper.

\bt{t-sigmao}
For any $\tau \ge 0$ and any $\alpha \in\R$, 
there is a generic subset $\Ss$ of $\sysI$
such that, for every system $\sigma ^o\in \Ss$, there is some
$\wh{\sigma ^o}\in\Ss$ such that 
\begin{enumerate}
\item $\sigma ^o$ and $\wh{\sigma ^o}$ are i/o equivalent under the
pulse function $\pulseu$ ($\sigma ^o \ioeq{\{\pulseu\}}\wh{\sigma ^o}$), but
\item $\sigma ^o$ and $\wh{\sigma ^o}$ are not i/o equivalent.
\end{enumerate}
\ets

\bt{t-sigma}
There is a generic subset ${\cal G}$ of $\sysII$
such that, for every system $\sigma \in {\cal G}$, there is some
$\wh{\sigma }\in\Gg$ such that 
\begin{enumerate}
\item $\sigma $ and $\wh{\sigma }$ are i/o equivalent under all the
pulses in the set $\Uu$ ($\sigma ^o \ioeq{\Uu}\wh{\sigma ^o}$), but
\item $\sigma $ and $\wh{\sigma }$ are not i/o equivalent.
\end{enumerate}
\ets

By setting $\tau =0$ for the collection $\Uu$, 
one obtains the following as a consequence of Theorem~\ref{t-sigma}: 

\bc{c-sigma} 
There is a generic subset ${\cal G}$ of $\sysII$
such that, for every system $\sigma \in {\cal G}$, there is some
$\wh{\sigma }\in\Gg$ such that 
\begin{enumerate}
\item $\sigma $ and $\wh{\sigma }$ 
are i/o equivalent under every constant input, but
\item $\sigma $ and $\wh{\sigma }$ are not i/o equivalent.~\mybox
\end{enumerate}
\ecs

The first part of Corollary \ref{c-sigma} may be restated as
follows: for every $(A, N, b, c)\in\Gg$, 
there exists some $(\wh A, \wh N, \wh b, \wh c)\in\Gg$
such that
\[
ce^{(A+\beta N)t}b \;=\; \wh ce^{(\wh A+\beta \wh N)t}\, \wh b\qquad
\]
for all $\beta \in \R$ and all $t\geq 0$.
This implies that:
\[
c \int_0^t e^{(A+\beta N)(t-s)}\beta b\,ds \;=\; 
\wh c\int_0^t e^{(\wh A+\beta \wh N)t}\, \beta \wh b\,ds
.
\]
for all $\beta \in \R$ and all $t\geq 0$.
Hence, the result in Corollary \ref{c-sigma} also applies to systems as
in \rref{e-syso}:
\bc{c-syso}
There is a generic subset ${\cal G}$ of $\sysI$
such that, for every $\sigma \in {\cal G}$, there is some
$\wh{\sigma }\in\Gg$ such that 
\begin{enumerate}
\item $\sigma $ and $\wh{\sigma }$
are i/o equivalent under every constant input, but
\item $\sigma $ and $\wh{\sigma }$ are not i/o equivalent.~\mybox
\end{enumerate}
\ecs

Next, we state our positive results for systems of both types.  For  any
$\alpha \in\R$, let $\Vv_{\alpha }$ denote the set of pulses of magnitude $\alpha $:
\[
\Vv_{\alpha } :=\{\pulseu|\; \tau \ge 0\}.
\]

\bt{t-p1} For each $\alpha \not=0$, there is a generic subset $\Mm$ of
$\sysI$ such that, for every pair of systems $\sigma ^o_1$, $\sigma ^o_2\in\Mm$,
\[
\sigma ^o \ioeq{\Vv_{\alpha }}\wh{\sigma ^o} \; 
\Longleftrightarrow \;
\sigma ^o \ioeQ \wh{\sigma ^o}\,.
\]
\ets

\bt{t-p2} For each $\alpha \not=0$, there is a generic subset $\Mm$ of
$\sysII$ such that, for every pair of systems $\sigma ^o_1$, $\sigma ^o_2\in\Mm$,
\[
\sigma ^o \ioeq{\Vv_{\alpha }}\wh{\sigma ^o} \; 
\Longleftrightarrow \;
\sigma ^o \ioeQ \wh{\sigma ^o}\,.
\]
\ets

\section{Proofs of Negative Results}
\label{proof-neg-results}

\subsection{Some Preliminaries}

The following construction is key to the proofs of the negative results.
The following observation was apparently first made in~\cite{brockett_book}
(see problem 1 in page 110, and problem 12 in page 105).

\bl{l-0}
For each canonical triple $\sigma =(A, b, c)$, there is a unique matrix
$T=T(\sigma )$ such that
\be{l-0-eqs}
AT=TA'\,, \; b=Tc'\,,\; cT=b' \,.
\ee
Moreover, the matrix $T(\sigma )$ is given by
$T(\sigma )= {\cal R}(A, b)\,[({\cal O}(A, c))']^{-1}$.
\els

\notacc{
\bpr 
Observe that for each canonical triple $(A,b,c)$, the triple $(A, c', b')$ is
also canonical, and the two triples are i/o equivalent since
$cA^kb=(cA^kb)'=b'(A')^kc'$ for all nonnegative integers $k$.
Thus there is a (unique) similarity between $(A,b,c)$ and $(A,c',b')$, an
invertible matrix $T$ such that:
\[ 
AT=TA'\,, \; b=Tc'\,,\; cT=b' \,.
\]
The formula for $T$ is given in~(\ref{lin-T}), which, since
in this case $\wh{A}=A'$ and $\wh{b}=c'$ and
${\cal O}(A, c)'= {\cal R}(A',c')$, reduces to that shown.
\epr

Although not needed, it is worth remarking that $T$ is symmetric.
This can be proved as follows: transposing the relations in~(\ref{l-0-eqs}),
one has that also
$AT'=T'A'$, $b=T'c'$, and $cT'=b'$.
Since the similarity $T$ is unique, $T=T'$.
}
For each nonzero $n\times n$ matrix $S$, consider the following 
\foracc{proper linear subspace of $\R^{n\times n}$:}
\notacc{set:
\[
}
\foracc{$}
{\cal B}(S):=\{N\in \R^{n\times n}\,|\, NS = SN'\}.
\foracc{$}
\notacc{\]}
\notacc{
Note that ${\cal B}(S)$ is a proper linear subspace of $\R^{n\times n}$,
because, in particular, when $E_{ij}$ is the 
matrix having a $1$ in its $(i,j)$th position and zero elsewhere,
$E_{ij}S=SE_{ij}'$ implies that $s_{kj} = 0$ for all $k\not = i$. 

}
We now define a set $\Gg_0$ that will play a major role in the constructions.
It is defined as the set consisting of those 4-tuples
\[
(A,N,b,c) \in  
\R^{n\times n}\times \R^{n\times n}\times \R^{n\times 1}\times \R^{1\times n}
\]
such that:
\begin{enumerate}
\item
$(A,b,c)$ is canonical,
\item
$N\not\in {\cal B}(T(A,b,c))$.
\end{enumerate}

Note that, since the triple $(A,b,c)$ is already canonical, every element of
$\Gg_0$ is canonical as a 4-tuple.

\notacc{%
For the next result, it is more elegant not to use inverses.
Given a triple $\sigma =(A,b,c)$, let $\tilde{\cal O}$ denote the cofactor matrix of
${\cal O}={\cal O}(A,c)$.  Note that if $\sigma $ is observable, then
${\cal O}^{-1} = \frac{1}{\Delta }\tilde{\cal O}$, where $\Delta (\sigma )=\det{\cal O}$.

For any triple $\sigma =(A,b,c)$ (not necessarily minimal), we define
$\tilde T(\sigma ):= {\cal R}(A, b)\,(\hat {\cal O}(A, c))'$.
This is a polynomial expression on the entries of $A$, $b$, and $c$.
If $(A,c)$ is observable, $T(\sigma )=\frac{1}{\Delta (\sigma )}\tilde T(\sigma )$.

Observe that ${\cal B}(S)={\cal B}(\delta S)$ for any scalar $\delta $.
Thus, ${\cal B}(T(A,b,c))={\cal B}(\tilde T(A,b,c))$.
}

\bl{l-go}
The complement $\Gg_0^c$ of $\Gg_0$ is a proper algebraic subset of 
$\R^{n\times n}\times \R^{n\times n}\times \R^{n\times 1}\times 
\R^{1\times n}$. 
\els

\notacc{%
\bpr
The complement of $\Gg_0^c$ of $\Gg_0$ is the union of the solution
sets of the following equations respectively:
\beqn
\det\,{\cal R}(A, b)&=& 0,\\
\det\,{\cal O}(A, c)&=& 0,
\eeqn
and the $n^2$ scalar equations given by
\be{e-tn}
N\tilde T(\sigma ) = \tilde T(\sigma )N'.
\ee
Hence, $\Gg_0^c$ is an algebraic set.
Each subset is proper (for the last one, pick an arbitrary canonical $(A,b,c)$
and refer to the above remark that ${\cal B}(S)$ is always proper), and
hence of dimension less than $n^2+2n$, so the union is also proper.
\epr
}

\bl{l-linear}
For each $(A,N,b,c)\in \Gg_0$, consider $M:=TN'T^{-1}$, where $T=T(A,b,c)$.
Then, 
\begin{enumerate}
\item $M\not= N$,
\item $(A, M, b, c)\in\Gg_0$, and
\item  for each $\gamma \in \R$ 
and each nonnegative integer $k$:
\be{e-constant}
c(A+\gamma N)^kb \;=\; c(A+\gamma M)^kb \,.
\ee
\end{enumerate}
\els

\notacc{
\bpr
Let  $(A,N,b,c)\in \Gg_0$.
The fact that $M\not=N$ follows from the fact that $N\not\in{\cal
  B}(T(A,b,c))$.  To see that $(A, M, b, c)\in\Gg_0$, note that
$MT = TN'T^{-1}T=TN' \not= TM'$ since $N\not= M$ and $T$ is invertible. 
The equality \rref{e-constant} follows by the equalities
$AT=TA'$, $b=Tc'$, $cT=b'$, $TN'=MT$, and the following:
\beqn
c(A+\gamma N)^kb 
&=& (c(A+\gamma N)^kb)' 
= b'(A'+\gamma N')^kc'\\
&=& cT\left(T^{-1}AT + \gamma T^{-1}MT\right)^kT^{-1}b
= c(A+\gamma M)^kb 
\eeqn
for all $k\ge 0$ and all $\gamma \in\R$.
\epr
}

\bc{so-not-equiv}
For each $(A,N,b,c)\in \Gg_0$, and $M=TN'T^{-1}$, the 4-tuples
$(A, N, b, c)$ and
$(A, M, b, c)$ and
are not i/o equivalent to each other.
\ecs

\bpr
Suppose that these two 4-tuples would be i/o equivalent.
By Lemma~\ref{main-theo-canonical}, they are similar.
Let $T$ provide a similarity as in~(\ref{e-equiv}).
In particular, $T$ provides a similarity between the canonical triple
$(A,b,c)$ and itself.  Since there is a unique such similarity, and
the identity $I$ is one, it follows that $T=I$.  Thus 
$N = TMT^{-1}=M$, contradicting the fact that $N\not= M$.
\epr

\subsection{Proof of \protect{Theorem~\ref{t-sigmao}}}

Let $\Cc$ be the subset consisting of all those 4-tuples
\[
(Q, N, b_0, c)\;\in \;\R^{n\times n}\times \R^{n\times n}\times \R^{n\times 1}\times \R^{1\times n}
\] 
which satisfy the following conditions:
\begin{enumerate}
\item[(a)] $(Q, b_0, c)$ is canonical,
\item[(b)] $(Q-N, b_0, c)$ is canonical, 
\item[(c)] $e^Q - I$ is invertible, and
\item[(d)] $N\notin {\cal B}(T(Q,b_0,c))$.
\end{enumerate}

Letting $\Lambda (F)$ denote the collection of eigenvalues of a matrix $F$,
the Spectral Mapping Theorem implies that
$e^{\Lambda (Q)} = \{e^\lambda : \; \lambda \in \Lambda (Q)\}$.
Thus, assumption (c), which says that $1$ is not an eigenvalue of $e^Q$,
implies, in particular, that $Q$ is invertible.

\bl{l-C} 
The complement $\Cc^c$ of $\Cc$ is a countable union
$\Cc^c=\bigcup _{k=0}^\infty \Ee_k$, 
where each $\Ee_k$ is a proper
algebraic subset of $\R^{n\times n}\times \R^{n\times n}\times 
\R^{n\times 1}\times \R^{1\times n}$.  
\els

\bpr
First note that $\Cc\subseteq \Gg_0$, and
$\Cc^c$ is the union of $\Gg_0^c$ and of the solution sets of
of the following equations:
\beq
& &\det{\cal R}(Q-N, b_0) = 0\label{e-C21},\\
& &\det{\cal O}(Q-N, c) = 0\label{e-C22},\\
& &\det\,(e^Q - I) = 0. \label{e-C4}
\eeq
Clearly, the solutions sets $\Aa_1$ and $\Aa_2$ of
Equations \rref{e-C21} and \rref{e-C22} 
respectively
are proper algebraic sets.  By the Spectral
Mapping Theorem,  \rref{e-C4} holds if and only if 
$-4k^2\pi ^2$ is an eigenvalue of $Q^2$ for some integer $k$, and hence,
the solution set $\Aa_3$ of 
\rref{e-C4} is the countable union of the solution sets $\{\Aa_{3k}\}$ 
of the equations
\[
\det(Q^2 + 4k^2\pi ^2I) = 0, \qquad k=0, 1, 2,\ldots .
\]
Hence, $\Cc^c$ is the countable union of  $\Gg_0^c$, $\Aa_1$, $\Aa_2$, and
$\{\Aa_{3k}\}_{k\ge 0}$. 
\epr

Let $\Xx = \R^{n\times n}\times \R^{n\times n}\times 
\R^{n\times  1}\times \R^{1\times n}$, and
consider the analytic map $\psi : \Xx\rightarrow \Xx$
defined by
\[
\psi : (Q,N,b_0,c) \mapsto  \left(Q-N, \, N, \, 
(\rho (Q))^*b_0, c\right),
\]
where $\rho (Q) = \int_0^1e^{sQ}ds$, and $\rho (Q)^*$ denotes 
the adjoint matrix of $\rho (Q)$.
Note that if $e^Q - I$ is invertible, then $Q$ is
invertible, and the matrix $\rho (Q) = Q^{-1}(e^Q-I)$ is also
invertible.  Hence, when restricted to
the open set
\[
\Xx_1 \,:=\; \left\{(Q, N, b_0, c):\;\mbox{$e^Q-I$ is invertible}\right\}\,,
\]
$\psi $ is given by
\be{e-psio}
 (Q,N,b_0,c) \mapsto  \left(Q-N, \, N, \;\det(\rho (Q))
\left(\rho (Q)\right)^{-1}\!
b_0, \; c\right).
\ee
Let $\psi _0$ denote the restriction of $\psi $ to $\Xx_1$, and
consider the open set
\[
\Yy_1 \,:=\;
  \left\{(A, N, b, c):\;\mbox{$(e^{A+N}-I)$ is invertible}\right\} \,.
\]
Then, $\psi _0$ is a (smooth)  diffeomorphism from $\Xx_1$ to 
$\Yy_1$.
\notacc{%
Its inverse is $\psi _0^{-1}: \Yy_1\rightarrow \Xx_1$ given by
\be{e-psio-inverse}
\psi _0^{-1}: \; (A, N, b, c) \mapsto 
\left(A+N, \, N, \, 
\left[\det\!\left(\int_0^1e^{s(A+N)}ds\right)\right]^{-1}\!
\!\left(\int_0^1e^{s(A+N)}ds\right)
b, \, c\right).
\ee
}
Since $\Yy=\psi (\Xx_1)\subseteq \psi (\Xx)$, it follows that the complement
$(\psi (\Xx))^c$ of $\psi (\Xx)$ is a subset of 
$\{(A, N, b, c):\;\mbox{$(e^{A+N}-I)$ is singular}\}$.
Hence, $(\psi (\Xx))^c$ is contained in the countable union of the solution
sets $\Ff_k$ of the equations
\[
\det((A+N)^2 + 4k^2\pi ^2I)=0, \qquad k=0, 1, 2,\ldots \,.
\]
Let $\Dd= \psi (\Cc)$, and
write $\Cc^c = \bigcup _{k=0}^\infty \Ee_k$, where the sets $\Ee_k$ are as in
Lemma~\ref{l-C}.
The next lemma then follows from the fact
that
$\Dd^c\subseteq [\psi (\Xx)]^c\bigcup \psi (\Cc^c)$, and Lemma~\ref{l-C}:

\bl{l-cD}
$\Dd^c \subseteq\Bigl(\bigcup _{k=0}^\infty\psi _0(\Ee_k)\Bigr)\bigcup 
 \Bigl(\bigcup _{k=0}^\infty \Ff_k\Bigr)$, where $\psi _0$ is the
diffeomorphism from $\Xx_1$ to $\Yy_1$ defined by \rref{e-psio}.
\el

\bc{c-D} 
$\Dd$ is generic.
\ec

\bpr
Since every proper algebraic set has measure zero, the set
$\bigcup _{k\ge 1}\Ff_k$ has measure zero.  Furthermore, since the image of 
a measure zero set under a differentiable map has measure zero (see
e.g.~\cite{Brocker}, Lemma 2.6), 
$\psi (\Ee_k)$ has measure zero for each $k$.  This implies
that $\Dd^c$ has measure zero.  Therefore, $\Dd$ is of full measure, and
as a consequence, $\Dd$ is dense.
Finally, $\Dd$ is open because it is the
image of $\Cc$ under the diffeomorphism $\psi _0$
and $\Yy_1$ is an open subset of
$\R^{n\times n}\times \R^{n\times n}\times \R^{n\times 1}\times \R^{1\times n}$.
\epr

Let ${\bf u} = \pulseu$ with $\tau =1, \alpha =1$.
\bl{l-syso}
Consider systems as in \rref{e-syso}.
For every  $\sigma ^o\in\Dd$,
there exists $\wh{\sigma ^o}\in \Dd$ such that the following holds:
\begin{enumerate}
\item $\sigma ^o$ and $\wh{\sigma ^o}$ are i/o equivalent under
  the pulse function ${\bf u}$, but
\item $\sigma ^o$ and $\wh{\sigma ^o}$ are not i/o equivalent.
\end{enumerate}
\els

\bpr
Let $\sigma ^o=(A, N, b, c)\in\Dd=\psi ({\cal C})$.  Then there exists $(Q, N, b_0, c)\in
\Cc\subseteq\Gg_0$ such that 
\[
A = Q-N, \quad
b=[\det(\rho (Q))](\rho (Q))^{-1}b_0.
\]
Let $b_1 = \rho (Q)b$.  Then
$b_1 = \det(\rho (Q))b_0$, and hence,
\beqn
{\cal R}(Q, b_1) &=& \det(\rho (Q)) {\cal R}(Q, b_0)\\
{\cal R}(Q-N, b_1) &=& \det(\rho (Q)) {\cal R}(Q-N, b_0).
\eeqn
Since $\det(\rho (Q)) \not=0$, both $(Q, b_1)$ and $(Q-N, b_1)$ are
reachable.  Moreover,
\beqn
T(Q, b_1, c) &=& T(Q, \det(\rho (Q))b_0, c)\\
&=& {\cal R}(Q, \det(\rho (Q))b_0) [{\cal O}(Q, c)']^{-1}\\
&=&
\det(\rho (Q))\,{\cal R}(Q, b_0)[{\cal O}(Q, c)']^{-1}\\
&=&
 \det(\rho (Q))\,T(Q, b_0, c),
\eeqn
which implies that
${\cal B}(T(Q, b_1, c)) = {\cal B}(T(Q, b_0, c))$.
Therefore, $(Q, N, b_1, c)\in\Cc$.
In particular, $(Q,N,b,c)\in \Gg_0$.

Applying Lemma  
\ref{l-linear} to $(Q, N, b_1, c)$, one sees that with $M =TN'T^{-1}$
(where $T = T(Q, b_1, c)$), it holds that  $M\not=N$, and 
\[
ce^{t(Q+\gamma N)}b_1 = ce^{t(Q+\gamma M)}b_1 \;\;
\forall\, \gamma \,,\;\forall \,t\geq 0.
\]
In particular, for $\gamma =-1$,
\[
ce^{t(Q-N)}b_1 = ce^{t(Q-M)}b_1 \qquad \forall \,t\geq 0\,.
\]
Let $\wh{\sigma ^o} = ((A+N-M), M, b, c)$.  Consider the two systems
$\sigma ^o$ and  $\wh{\sigma ^o}$:
\beqn
\dot x &=& (A + uN)x + bu, \;\; x(0)=0, \;\; y = cx,\\
\dot z &=& [(A+N-M) + uM]z + bu, \;\; z(0) = 0, \;\; y = cz.
\eeqn
Since ${\bf u}=1$ for $t\in[0, 1]$, the two systems reduce in that interval to:
\beqn
\dot x &=& (A+ N)x + bu, \quad x(0)=0, \;\;  y = cx, \;\;  0\le t\le 1\\
\dot z &=& (A+N)z + bu, \quad z(0)=0,\;\;  y=cz, \;\;  0\le t\le 1.
\eeqn
It follows that $x(t)=z(t)$ on $[0, 1]$, and hence outputs coincide for
$t\in [0,1]$.
In particular, at time $t=1$ both systems are in state
\[
\int_0^1 e^{(1-s)(A+N)}b\,ds = \int_0^1 e^{(1-s)Q}b\,ds
 = \rho (Q)b
= b_1 \,.
\]
Now, for $t\geq 1$, using that $A=Q-N$ and $A+N-M=Q-M$, we have that
\[
\dot x = (Q-N)x\,,\quad x(1) = b_1\,,
\]
\[
\dot z = (Q-M)z\,,\quad z(1) = b_1\,,
\]
and, since $ce^{t(Q-N)}b_1 \equiv  ce^{t(Q-M)}b_1$, it follows that the
outputs are the same for all $t>1$ as well.

To show that $\sigma ^o$ and $\wh{\sigma ^o}$ are not equivalent, we will
show explicitly that the inputs 
\[
u_s(t) = \left\{
  \begin{array}{ll}
1, & 0\le t\le 1,\\
0,& 1 < s\le 1+s,\\
1, & t> 1+s
  \end{array}
\right.
\]
(with varying $s\ge 0$) are enough to distinguish the two systems.

Suppose that the systems have the same output functions under these
input functions.

Let $s\ge 0$ and $t\ge0$.
Since $x(1)=z(1)=b_1$, and $x(1+s)=e^{s(Q-N)}b_1$,
$z(1+s)=e^{s(Q-M)}b_1$, at time $t> 1+s$ we have:
\[
\dot x = Qx + b\,,\; x(1+s)=e^{s(Q-N)}b_1,
\]
\[
\dot z = Qz + b\,,\; z(1+s)=e^{s(Q-M)}b_1 \,.
\]
{}From $cx(t)\equiv cz(t)$ and taking derivatives, one gets
\[
cQx(t) +cb= cQz(t) + cb\qquad\forall\,t>1+s.
\]
Hence,
\[
cQx(t) =  cQz(t)\qquad\forall\,t>1+s\,.
\]
Taking more derivatives, we conclude inductively that
\[
cQ^kx(t) =  cQ^kz(t)\qquad\forall\,t>1+s, \ \forall\,k\ge 0.
\]
In particular by continuity, we have that
$cQ^kx(1+s)=cQ^kz(1+s)$ for all $k\ge 0$, and by 
observability of $(Q, c)$, 
\[
x(1+s)=z(1+s) \qquad\forall\,s\ge 0,
\]
that is,
\[
e^{s(Q-N)}b_1 = e^{s(Q-M)}b_1\qquad\forall\, s\ge 0.
\]
Taking derivatives with respect to $s$, we conclude that:
\[
(Q-N)^kb_1 = (Q-M)^kb_1 \qquad\forall\,k\ge 0.
\]
This implies that
\[
(Q-N){R}_1 = (Q-M){R}_2,\quad{\rm and}\ R_1 = R_2,
\]
where ${R}_1 = {\cal R}(Q-N, b_1)$, and ${R}_2 = {\cal R}(Q-M,
b_1)$.  Since $(Q-N,b_1)$ is reachable, 
$R_1$ is invertible. From this we conclude that $Q-N=Q-M$, and hence
$N= M$, a contradiction. 

To complete the proof of Lemma \ref{l-syso}, we show that
$\wh{\sigma ^o}\in \Dd$.  First observe that:
\beqn
\wh{\sigma ^o} &=& (A+N-M, M, b, c) \;=\; (Q-M, M, b, c)\\
&=&\psi _0(Q, M, b_0, c).
\eeqn
Thus, if we prove that $(Q, M, b_0, c)\in\Cc$,
then $\wh{\sigma ^o}=\psi (Q, M, b_0,c) \in\Dd$.
By Lemma
\ref{l-linear}, $(Q, M, b_0, c)\in\Gg_0$.  It is thus enough to show
that $(Q-M, b_0, c)$ is canonical.  
To see this, note that since $b_1 = \det(\rho (Q))b_0$ and $((Q-N), b_0, c)$ is
canonical, it follows that $(Q-N, b_1, c)$ is canonical. Thus, $(Q-M,
b_1, c)$ is canonical as $(Q-M, b_1, c)$ is similar to $((Q-N)', c',
b_1')$.  Again, applying the fact that  $b_1=
\det(\rho (Q))b_0$, one sees that $(Q-M, b_0, c)$ is canonical.
\epr

The above completes the proof of Theorem~\ref{t-sigmao} for the special case
$\tau =1$ and $\alpha =1$.  The general case can be obtained by rescaling inputs and
time scale, as follows.
We consider the impulse function $\pulseu$ for any fixed $\tau >0$ and
$\alpha \in\R$.  Without loss of generality, we assume that $\alpha \not=0$.
For the initial-value problem 
\[
\dot x = (A+\alpha N)x + \alpha b\,, \;\;x(0) = 0,
\]
let $N_\alpha = \alpha N, b_\alpha = \alpha b$, and 
consider the initial-value problem
\[
\dot {\tilde x} = \frac{1}{\tau }\Bigl((A + N_\alpha )\tilde x +  b_\alpha \Bigr), \;\;
\tilde x(0) = 0.
\]
Then $\tilde x(t) = x(t/\tau )$.
It then can be seen that, with
\[
\Dd_{\alpha , \tau }:=
\left\{\left(\frac{1}{\tau }A,\, \frac{\alpha }{\tau }N, \,
\frac{\alpha }{\tau }b, \, c\right):\; (A, N, b, c)\in\Dd\right\}
\]
and any $\sigma \in\Dd_{\alpha , \tau }$, there exists some $\wh
\sigma \in\Dd_{\alpha , \tau }$ such that $\sigma $ and $\wh\sigma $ have
the same output under the impulse input $\pulseu$, but the two systems are not
equivalent. 

\subsection{Proof of \protect{Theorem~\ref{t-sigma}}}

In this section we consider systems defined as in \rref{e-sys}.
Let $\tau \ge 0$ and $\alpha \in\R$ be given.
Consider the analytic map
\notacc{%
\[
\Phi : \R^{n\times n}\times \R^{n\times n}\times \R^{n\times 
  1}\times \R^{1\times n}
\rightarrow \R^{n\times n}\times \R^{n\times n}\times \R^{n\times 
  1}\times \R^{1\times n}
\]
given by
}
\[
\foracc{\Phi : \;\;}
(P, \, N, \, b_0, \,c) \mapsto (P - \alpha N, \, N, \, e^{-\tau P}b_0, \, c).
\]
This is an analytic diffeomorphism whose inverse map is given by
\[
(A, N, b, c)\mapsto (A + \alpha N, \, N, \, e^{\tau (A+\alpha N)}b, \, c).
\]
Let $\Hh = \Phi (\Gg_0)$.  The following is a consequence of Lemma \ref{l-go}:

The complement $\Hh^c$ of $\Hh$ is the image of an proper algebraic set under
a diffeomorphism from $\Xx$ to $\Xx$: $\Hh^c = \Phi (\Gg_0^c)$.
Since the image of a measure zero set under a smooth map 
(see e.g.~\cite{Brocker}, Lemma 2.6) has measure zero, we conclude:

\bc{c-h}
The collection $\Hh$ is generic.
\ec

\bl{l-sigma}
For any $\sigma _1\in\Hh$, there exists $\sigma _2\in\Hh$ such that
\begin{enumerate}
\item[(a)] $\sigma _1$ and $\sigma _2$ have the same output for any $u\in
  \Uu$; 
\item[(b)] the two systems $\sigma _1$ and $\sigma _2$ are not equivalent.
\end{enumerate}
\els

Let $\sigma _1=(A,N,b,c)\in\Hh$.
Thus, by definition of $\Hh$, there exists $(P, N, b_0, c)\in\Gg_0$ such
that $\sigma _1 = (P - \alpha N, N, e^{-\tau P}b_0, c)$.  
Let $M =TN'T^{-1}$, where $T = T(P, b_0, c)$.  Then $M\not= N$,
$(P, M, b_0, c)\in\Gg_0$, and
the system $\sigma _2:=(P-\alpha M, M,  e^{-\tau P}b_0, c) =\Phi (P, M,
b_0, c) \in\Hh$.

To prove part (a), pick any $u\in {\cal I}_{\tau , \alpha }$, and assume
$u(t) = \beta $ for $t > \tau $.  The two systems are given by
 \beqn
 \dot x &=& ((P - \alpha N) \,+\, uN)x, \quad x(0) =b,\\
 \dot z &=& ((P - \alpha M) + uM)z, \quad z(0) = b.
 \eeqn
For $0\le t\le \tau $, both systems reduce to the same equation:
\[
\dot p = P \, p, \quad p(0) = b
= e^{-\tau P}b_0,
\]
and in particular, $cx(t) = cz(t) = ce^{tP}b$ for all $0\le t\le \tau $.
For $t > \tau $, the two systems become
\beqn
\dot x &=& ((P-\alpha N) +\, \beta N)x, \quad x(\tau ) = 
b_0\\
\dot z &=& ((P-\alpha M) + \beta M)z, \quad z(\tau )= 
b_0
\eeqn
and thus have the respective solutions
\beqn
x(t) &=& ce^{(P +(\beta - \alpha )N)(t - \tau )}x(\tau )
= ce^{(P +(\beta - \alpha )N)(t - \tau )}b_0,\\
z(t) &=& 
 ce^{(P +(\beta - \alpha )M)(t - \tau )}z(\tau )
= ce^{(P +(\beta - \alpha )M)(t - \tau )}b_0.
\eeqn
By Lemma \ref{l-linear} and the choices of $M$ and $N$, it follows with
$\gamma = \beta -\alpha $ that 
\[
ce^{(P + \gamma N)(t - \tau )}b_0
=
ce^{(P + \gamma M)(t - \tau )}b_0,
\]
which implies that $cx(t) = cz(t)$ for all $t\ge \tau $ as well, so part (a)
is proved.

To prove part (b), suppose the two systems $\sigma _1$ and
$\sigma _2$ have the same output for all inputs, and so, in particular, for all
inputs $u$ for which $u(t) = \alpha $ for $0\le t\le \tau $.  This implies
that the two systems 
\beqn
\dot x &=& ((P-\alpha N) \,+\, u N)x,\quad x(\tau ) = e^{\tau P}b,\quad
y = cx,\\
\dot z &=& ((P-\alpha M) + u M)z, \quad z(\tau ) = e^{\tau P}b,\quad 
y=cz,
\eeqn
have the same
output for any input.
Rewriting  the two systems as:
\beqn
\dot x &=& (P\,+\, (u-\alpha ) N)x,\quad x(\tau ) = e^{\tau P}b,\quad
y = cx,\\
\dot z &=& (P \, + (u-\alpha ) M)z, \quad z(\tau ) = e^{\tau P}b,\quad 
y=cz,
\eeqn
and writing $v(t) = u(t)- \alpha $ for a new input
$v$, one sees that the two systems
\beqn
\dot x &=& (P + vN)x, \quad x(\tau ) = b_0, \quad y = cx,\\
\dot z &=& (P + vM)z, \quad z (\tau ) = b_0, \quad y = cz,
\eeqn
have the same outputs for all inputs $v$ and times $\geq \tau $, which is the same
as saying that the two systems of type II with associated 4-tuples
$(P,N.b_0,c)$ and $(P,M,b_0,c)$ are i/o equivalent,
which is a contradiction in view of Lemma~\ref{l-linear}
and Corollary~\ref{so-not-equiv}.

\notacc{ 
This completes the proof of Theorem~\ref{t-sigma}.
}

\section{Proofs of Positive Results}
\label{proof-pos-results}

In this section we prove the positive results .

Let $\alpha \not=0$ be given.  Let $\Mm$ be the set of 4-tuples
satisfying the following two properties:
\begin{enumerate}
\item $(A, b, c)$ is canonical,
\item $(A+\alpha N, b)$ is controllable.
\end{enumerate}
Since the complement $\Mm^c$ of $\Mm$ is defined by the union of
the solution sets of the equations
\beqn
\det\,{\cal R}(A ,b)&=& 0,\\
\det\,{\cal O}(A, c)&=& 0,\\
\det\,{\cal R}(A+\alpha N, \, b)&=&0,
\eeqn
we have:

\bl{l-old1}
The complement $\Mm^c$ of $\Mm$ is a proper algebraic set.
Consequently, $\Mm$ is open, dense, and of full measure.
\els

We will prove that the sets of systems in $\sysI$ and $\sysII$ whose 4-tuples
are in $\Mm$ satisfy the conclusions of Theorems \ref{t-p1} and
\ref{t-p2} respectively.

\subsection{Proof of \protect{Theorem \ref{t-p1}}}

\bl{l-old2}
Assume that $(A, b)$ is controllable. Then,
for any $N\in\R^{n\times n}$ and any $\alpha \in\R$, 
\[
\left(A, \, \left(\int_0^\tau e^{(A+\alpha N)s}ds\right) b\right)
\]
is controllable for almost all $\tau >0$.
\els

\bpr Let $A, N , b$ be given.  For any $\tau >0$,
\[
\int_0^\tau e^{(A+\alpha N)s}\,ds
=\sum_{k=0}^\infty\frac{(A+\alpha N)^k}{(k+1)!}\tau ^{k+1}
= \tau \Psi (\tau ),
\]
where $\Psi (\tau ) = 
\sum_{k=0}^\infty\frac{(A+\alpha 
    N)^k}{(k+1)!}\tau ^{k}$.
Since $\Psi (0) = I$, $\Psi $ is analytic, and $(A, b)$ is controllable,
it follows that $(A, \Psi (\tau )b)$ is controllable for almost all
$\tau $.  It then follows that
$(A, \tau \Psi (\tau )b)$ is controllable for almost all $\tau >0$.
\epr

To prove Theorem \ref{t-p1}, we pick two
systems $\sigma _1^o =(A_1, N_1, b_1, c_1)$
and $\sigma _2^o=(A_2, N_2, b_2,c_2)$ in $\Mm$, and suppose that
they produce have the same
output function for each $u\in\Vv_\alpha $.
We must show that $\sigma _1^o\ioeQ\sigma _2^o$.

Fix a $\tau >0$.  Applying $\pulseu\in\Vv_\alpha $ to the two systems:
\beqn
\dot x &=& (A_1+u N_1)x + b_1 u,\quad
x(0)=0, \quad y=c_1x\\
\dot z &=&  (A_2+u N_2)z + b_2 u,\quad
z(0)=0, \quad y=c_2z,
\eeqn
one has:
\be{e-lem5.3.1}
c_1e^{A_1(t-\tau )}x(\tau ) = c_2e^{A_2(t-\tau )}z(\tau )\qquad\forall\,t\ge\tau ,
\ee
where
\beqn
x(\tau ) &=&
\notacc{\alpha \int_0^\tau e^{(A_1+\alpha N_1)(\tau -s)}ds\, b_1 =}
\alpha \int_0^\tau e^{(A_1+\alpha N_1)s}ds\, b_1,\\
z(\tau ) &=& 
\notacc{\alpha \int_0^\tau e^{(A_2+\alpha N_2)(\tau -s)}ds\, b_2 =}
\alpha \int_0^\tau e^{(A_2+\alpha N_2)s}ds\, b_2.
\eeqn
This holds for any $\tau >0$.

By Lemma \ref{l-old2}, we may pick some $\tau _0>0$ such that
$(A_1, x(\tau _0), c_1)$ and $(A_2, z(\tau _0), c_2)$ are both canonical.
Since by~(\ref{e-lem5.3.1}) these two triples are i/o equivalent,
there exists some invertible matrix $T\in\R^{n\times n}$ such that 
\be{e-aT}
A_2=T^{-1}A_1T, \;
z(\tau _0) = T^{-1}x(\tau _0),\;
c_2 = c_1T.
\ee
Using in particular that $c_2e^{A_2s} = c_1e^{A_1s}T$ for all $s$,
\rref{e-lem5.3.1} becomes:
\[
c_1e^{A_1(t-\tau )}x(\tau )
=c_1e^{A_1(t-\tau )}Tz(\tau )\qquad\forall\, t\ge\tau .
\]
{}From the observability of $(A_1, c_1)$, it follows that
\[
x(\tau ) = Tz(\tau )
\]
for all $\tau >0$.   Equivalently:
\[
\int_0^\tau e^{(A_1+\alpha N_1)s}ds\, b_1
=
T\int_0^\tau e^{(A_2+\alpha N_2)s}ds\, b_2
\]
for all $\tau >0$. Taking the derivative with respect to $\tau $, one gets:
\[
e^{(A_1+\alpha N_1)\tau } b_1 = T e^{(A_2+\alpha N_2)\tau }
b_2.
\]
Note that this is true for all $\tau \ge 0$.  
In particular,
\be{b-equals}
b_1 = Tb_2.
\ee
On the other hand,
taking repeated derivatives in
$\tau $ and then setting $\tau =0$, one obtains:
\be{e-a6}
(A_1+\alpha N_1)^kb_1 = T(A_2+\alpha N_2)^kb_2\qquad\forall\,k\ge 0.
\ee
This implies, with $0\le k\le n-1$,
\be{e-a4}
{\cal R}(A_1+\alpha N_1, b_1) = T[{\cal R}(A_2+\alpha N_2, b_2)],
\ee
and with $1\le k\le n$,
\beqn
(A_1+\alpha N_1)[{\cal R}(A_1+\alpha N_1, b_1)]
&=&\\
&&
\hskip-2cm
T(A_2+\alpha N_2)[{\cal R}(A_2+\alpha N_2, b_2)].
\eeqn
Combining this with \rref{e-a4}, one sees that
\beqn
(A_1+\alpha N_1)[{\cal R}(A_1+\alpha N_1, b_1)]
&=&\\
&&
\hskip-3cm
T(A_2+\alpha N_2)T^{-1}[{\cal R}(A_1+\alpha N_1, b_1)]
\eeqn
It then follows from the fact that ${\cal R}(A_1+\alpha N_1, b_1)$ is
invertible (because $(A_1+\alpha N_1, b_1)$ is controllable) that 
\[
T(A_2+\alpha N_2)T^{-1} = (A_1+\alpha N_1).
\]
It then again follows from \rref{e-aT} and the fact that $\alpha \not=0$
that $N_2 = T^{-1}N_1T$.  Combined with~(\ref{e-aT}) and~(\ref{b-equals}),
we have that the systems $\sigma _1^o$ and
$\sigma _2^o$ are similar, with the similarity matrix given by $T$, and
this completes the proof of Theorem~\ref{t-p1}.

\subsection{Proof of \protect{Theorem \ref{t-p2}}}

The proof of Theorem~\ref{t-p2} is almost the same as that of
Theorem~\ref{t-p1}, with Lemma \ref{l-old2} modified to the following: 

\bl{l-old4}
Assume that $(A, b)$ is controllable. Then,
for any $N\in\R^{n\times n}$ and any $\alpha \in\R$, 
\[
\left(A, \, e^{(A+\alpha  N)\tau } b\right)
\]
is controllable for almost all $\tau >0$.
\el

We pick elements $\sigma _1 =(A_1, N_1, b_1, c_1)$ and $\sigma _2=(A_2, N_2, b_2,c_2)$
of $\Mm$ (seen as a class of systems of type II)
which have the same output function for each $u\in\Vv_\alpha $. 

For any $\tau >0$, applying $\pulseu\in\Vv_\alpha $ to the two systems:
\beqn
\dot x &=& (A_1+u N_1)x,\quad
x(0)=b_1, \quad y=c_1x\\
\dot z &=&  (A_2+u N_2)z, \quad
z(0)=b_2, \quad y=c_2z,
\eeqn
one has, for any $\tau >0$,
\be{e-lem5.3.2}
c_1e^{A_1(t-\tau )}x(\tau ) = c_2e^{A_2(t-\tau )}z(\tau )\qquad\forall\,t\ge\tau ,
\ee
where
\[
x(\tau ) =  e^{(A_1+\alpha N_1)\tau }b_1,\;
z(\tau ) =  e^{(A_2+\alpha N_2)\tau } b_2.
\]
By Lemma \ref{l-old4}, there exists some $\tau _0>0$ such that both
$(A_1, x(\tau _0), c_1)$ and $(A_2, z(\tau _0), c_2)$ are canonical.
So, there exists some invertible matrix $T\in\R^{n\times n}$ such that 
\[
A_2=T^{-1}A_1T, \;
z(\tau _0) = T^{-1}x(\tau _0),\;
c_2 = c_1T,
\]
and consequently, \rref{e-lem5.3.2} becomes
\[
c_1e^{A_1(t-\tau )}x(\tau )
=c_1e^{A_1(t-\tau )}Tz(\tau )\qquad\forall\, t\ge\tau .
\]
For each $\tau >0$ given, using the observability of $(A_1, c_1)$, one
sees that :
\[
e^{(A_1+\alpha N_1)\tau }b_1
 =  
T
e^{(A_2+\alpha N_2)\tau } b_2
\]
for all $\tau >0$. Starting from here, one can complete the proof by
following the same steps as in the proof of Theorem \ref{t-p1}.

\notacc{
\section{A Remark on Sampled Controls}
\label{sec-juang}

As remarked earlier, our proofs of the positive results, Theorems
\ref{t-p1} and \ref{t-p2}, were inspired by the identification
algorithm presented in \cite{bilinear}.
That algorithm aims to find a system equivalent to the system being
identified, on the basis of observations at discrete instants
$0, \tau , 2\tau , \ldots , j\tau , \ldots $, where $\tau $ is a \emph{fixed} sampling time,
and having applied inputs which have the form $\pu{k\tau }$ (for varying
nonnegative integer $k$'s), i.e., pulses of magnitude $\alpha $ whose width is a
multiple of this same sampling time $\tau $.
The motivation is clear: one wishes to use a sample-and-hold strategy,
which is especially convenient for computer algorithms.
Unfortunately, this restriction to fixed sampling times means that
the algorithm cannot work for generic classes of systems, as we show
here by means of a counterexample.
(Mathematically, the difficulty is that some of the steps of the algorithm
given in~\cite{bilinear} involve taking logarithms of matrices,
which is an ambiguous procedure, as the author himself points out on the
paper.) 
To show this shortcoming, for any given $\alpha >0$, we produce an open
class $\Bb_\alpha $ of 2-dimensional systems of type I (it is easy to generalize to
larger dimensions and to systems of type II) with the following properties:
for every system in $\Bb_\alpha $, there 
is some other system, which is not i/o equivalent to the original one, yet
cannot be distinguished by applying steps of magnitude $\alpha $ and sampled
in the above way (with fixed $\tau $).
Thus, our approach, in which $\tau $ is varied, is actually necessary.

For any system $\sigmao = (A, N, b, c)\in\sysI$, we denote by $\sigmao_\tau $
the discrete time system which results from sampling the system with $\tau $ as
the length of the sampling interval, and using input functions $u$ that are
constant over each sampling interval:
\be{e-smp0}
x_{k+1} = F_kx_k + u_kg_k, \quad x(0)=0, \quad
y_k = cx_k,
\ee
where $x_k = x(k\tau ), y_k = y(k\tau )$, $u_k$ is the value of $u$ over
the interval $(k\tau , (k+1)\tau )$ and
\beqn
F_k &=& e^{(A+ u_kN)\tau },\\
g_k &=& 
\int_0^\tau e^{(A+u_kN)(\tau -s)}b\,ds.
\eeqn
For disrete time systems as in \rref{e-smp0}, the i/o equivalence under
a collection of inputs is defined in the same manner as in the
continuous time case.

Let $\Bb_\alpha $ be the subset consisting of systems of type I for which
the 4-tuples 
\[
 (A, N, b, c)\in \R^{2\times 2}\times \R^{2\times 2}\times \R^{2\times 
  1}\times \R^{1\times 2}
\]
satisfy:
\begin{enumerate}
\item $(A+\alpha N, b, c)$ is canonical; and
\item $A + \alpha N$ has a pair of conjugate complex eigenvalues 
$r\pm si$, with  $s\not= 0$. 
\end{enumerate}
Since the set of 4-tuples for which the triple $(A+\alpha N, b, c)$ is canonical
is generic, and the set of 4-tuples
$(A, N, b, c)$ for 
which $A+\alpha N$ has a pair of nonzero conjugate complex eigenvalues
has a nonempty interior (because of continuity of eigenvalues on
matrix entries), the set $\Bb_\alpha $ contains an open set.

\bp{p-n1} For any $\tau >0$, and
for any $\sigmao\in \Bb_\alpha $, there exists some $\wh{\sigmao}\in\Bb_\alpha $ 
that
\begin{enumerate}
\item $\sigmao_\tau $ and $\wh{\sigmao}_\tau $ are i/o equivalent under 
the collection $\{\pu{k\tau }\}_{k\ge 0}$, but
\item $\sigmao$ and $\wh{\sigmao}$ are not i/o equivalent.
\end{enumerate}
\eps

\bpr Let $\tau >0$ be given, and consider
$\sigmao = (A, N, b, c)\in\Bb_\alpha $. Without loss of generality,
assume that $A+\alpha N$ is already in ``real Jordan canonical form'':
\[
A + \alpha N = \pmatrix{r&-s\cr s&r}
\]
for some 
$s\not=0$.
(If this were not the case, one may simply apply a similarity, and at the end
of the argument transform back to original coordinates.)
Let
\[
\Lambda _0 =\pmatrix{0&-2\pi /\tau \cr 2\pi /\tau &0}.
\] 
Choose an integer $l\not=0$ such that
\begin{itemize}
\item 
 $s+ \frac{2l\pi }{\tau }\not=0$;
and
\item $((A+\alpha N) + l\Lambda _0, \; b, \, c)$
is canonical.
\end{itemize}
A generic integer works; note that (since $n=2$):
\beqn
\det {\cal O}(A+\alpha N + l\Lambda _0, c) &= &
\det {\cal O}(A+\alpha N, c) + \det {\cal O}(l\Lambda _0, c)\\
&=&\det {\cal O}(A+\alpha N, c) + l\det {\cal O}(\Lambda _0, c).
\eeqn
and similarly for reachability.
Let $M$ be given by
\[
 M = N + \frac{l}{\alpha } \Lambda _0.
\]
That is, $M$ is chosen so that
\[
A + \alpha M = A + \alpha N + l\Lambda _0 =
\pmatrix{r&-s - \frac{2l\pi }{\tau } \cr 
s + \frac{2l\pi }{\tau }&r}.
\]
Note that both $(A+\alpha N)$ and $(A+\alpha M)$ are
is invertible since $r+s i\not=0$ and $r + (s+ \frac{2l\pi }{\tau })i\not=0$. 
Let $\wh{\sigmao} = (A, M, \wh{b}, c)$, where 
$\wh{b} = (A+\alpha M)(A + \alpha N)^{-1}b$.

Next we show that the two sampled systems $\sigmao_\tau $ and
 $\wh{\sigmao}_\tau $ are i/o equivalent under the collection
$\{\pu{k\tau }\}_{k\ge 0}$.

Consider the input function $\pu{k\tau }$ for some $k\ge 0$. Clearly the
two systems have the same state trajectory and same output function if
$k=0$ (the input is constantly zero).
So assume $k\ge 1$.  For $j\le k-1$, the two sampled systems are
given respectively by:
\beqn
x_{j+1} &=& F x_j + \alpha g, \quad x_0 = 0, \quad y_j = cx_j,\\
z_{j+1} &=& \wh{F} z_j + \alpha \wh{g}, \quad
z_0 = 0, \quad y_j =cz_j,
\eeqn
where 
\beqn
F &=& e^{(A+\alpha N)\tau } = e^{r \tau }
 \pmatrix{\cos s\tau &-\sin s\tau \cr
  \sin s\tau &\cos s\tau } \\[2mm]
&=& 
e^{r \tau }
\pmatrix{\cos (s + \frac{2l\pi }{\tau })\tau &
 -\sin(s + \frac{2l\pi }{\tau })\tau \cr
  \sin (s + \frac{2l\pi }{\tau }) \tau &\cos(s + \frac{2l\pi }{\tau })\tau } 
= e^{(A+\alpha M)\tau } 
=
\wh{F},
\eeqn
and
\beqn
g &=&\int_0^\tau e^{(A+\alpha N)(\tau - \theta )} b\,d\theta 
= (e^{(A+\alpha N)\tau } - I)(A+\alpha N)^{-1}b
\\
&=& (e^{(A+ \alpha M)\tau }-I)(A+\alpha M)^{-1}\wh b\\
&=& \int_0^\tau e^{(A+\alpha M)(\tau -\theta )} \wh b\,d\theta 
=\wh{g}
\eeqn
(where we have used the fact that $\int_0^\tau e^{Q(\tau -\theta )}\,d
\theta =
(e^{Q\tau }-I)Q^{-1}$ for any invertible $Q$).
Hence, $x_{j} = z_{j}$ for all $0\le j\le  k$. 
In particular, outputs coincide at all sampling times $t=j\tau $, $j\leq k$.
Over the interval
$[k\tau , \infty)$, the two
systems $\sigmao$ and $\wh{\sigmao}$ are given by
\[
\dot w= Aw, \quad w(k\tau ) = x(k\tau ) = z(k\tau ),
\quad
y =cw.
\]
This in particular implies that $x_j = z_j$ for all $j\ge k$
and therefore outputs also coincide at sampling times $t=j\tau $, $j>k$
(as well as for any time $t>k$).
Thus,
 $\sigmao_\tau $ and $\wh{\sigmao}_\tau $ are i/o equivalent under
$\pu{k\tau }$.

To show that the two systems $\sigmao$ and $\wh{\sigmao}$ are not i/o
equivalent, consider the constant input $u\equiv \alpha $.  Assume the
output functions of the two systems are the same under this constant
$u$.  Then
\[
c\int_0^t e^{(A+\alpha N)(t-\theta )}b\,\alpha \,d\theta 
=
c\int_0^t e^{(A+\alpha M)(t-\theta )}\,\wh b\,\alpha \,d\theta \qquad\forall\,
t\ge 0,
\]
and equivalently,
\[
c\left(e^{(A+\alpha N)t} - I\right)(A+\alpha N)^{-1}b
=
c\left(e^{(A+\alpha M)t} - I\right)(A+\alpha M)^{-1}\wh b\qquad\forall\,t\ge 0.
\]
Since $(A+\alpha M)^{-1}\wh b = (A+\alpha N)^{-1}b$, one gets
\[
c\left(e^{(A+\alpha N)t} - I\right)(A+\alpha N)^{-1}b
=
c\left(e^{(A+\alpha M)t} - I\right)(A+\alpha N)^{-1}b\qquad\forall\,t\ge 0.
\]
This in turn implies that
\[
\int_0^t ce^{(A+\alpha N)(t-\theta )}b\,d\theta 
=
\int_0^t ce^{(A+\alpha M)(t-\theta )} b\,d\theta \qquad\forall\,
t\ge 0.
\]
Taking derivatives with respect to $t$ repeatedly, and then setting $t=0$, one
obtains: 
\[
c(A+\alpha N)^j b = c(A+\alpha M)^j b\qquad\forall\, j\ge 0.
\]
Since both $(A+\alpha N, b, c)$ and $(A+\alpha M, b, c)$ are canonical,
it follows that the two systems are similar, contradicting the fact
that $A+\alpha N$ and $A+\alpha M$ have different pairs of eigenvalues. 

Finally we show that $\wh{\sigmao}\in\Bb_\alpha $.  Note that $M$ was
chosen so that $(A+\alpha M, c)$ is observable.  It is left to show that
$(A + \alpha M, \wh{b})$ is controllable.  
Note that $P=A+\alpha N$ and $Q=A+\alpha M$ commute.
 This is because $Q=P+\frac{2l\pi }{\tau }J$, and $PJ=JP$,
where
\[
J=\pmatrix{0&-1\cr1&0}.
\]
Using that $PQ=QP$ and therefore $QP^{-1}=P^{-1}Q$:
\[
{\cal R}(Q,\wh b) 
 = (\wh b \;\; Q\wh b) 
= (QP^{-1}b \;\; QQP^{-1}b)
= (QP^{-1}b \;\; QP^{-1}Qb)
= (QP^{-1}){\cal R}(Q,b)
\]
Since $M$ was chosen so that $(Q, b)$ is controllable, it
follows that ${\cal R}(Q, \wh b)$ is non-singular, and hence, 
$(A+\alpha M, \wh b)$ is controllable, as claimed.
\epr

}

\section{Conclusions and Final Remarks}
\label{sec-final}

For bilinear systems, we showed that step inputs are not enough for
identification, nor do single pulses suffice, but that that the family of all
pulses (of a fixed amplitude but varying widths) do suffice.
We presented results for single-input single-output systems, since one can
obviously identify a multiple-input multiple-output system by considering each
pair of input and output channels separately, and hence the family of pulses
also works for the general case.

We emphasize that we dealt in this paper with ideal noise-free conditions, and
ignored stochastic aspects and noisy data, because the underlying theoretical
questions of what is ultimately achievable are easiest to understand in a
deterministic setting.  Tools such as those here have been used, however, in
the formulation of identification algorithms from noisy data, for bilinear
systems~\cite{verdult-verhaegen}.  Nor did we deal here with questions of
computational and sample complexity.  However, the methods used are quite
constructive and indeed have appeared in the same context in~\cite{bilinear},
where numerical implementations are studied; regarding sample complexity, we
leave for further research the generalization of learning-theory
results~\cite{dasgupta,kuusela} from the linear case to the classes of systems
considered here.

Finally, bilinear systems were picked because an elegant result can be
established for them, as well as their applicability and general interest.
However, the study of similar problems to those treated here for more general
classes of systems is of great interest.

\notacc{\newpage}


\end{document}